\documentclass[11pt]{article}
\usepackage{amsfonts}
\usepackage{latexsym,tikz}
\usepackage{amsmath}
 
\title{Multiple Dedekind Zeta Functions}
\author{Ivan Horozov}

\usepackage{color}
\date{}

\newtheorem{theorem}{Theorem} 
\newtheorem{lemma}[theorem]{Lemma}
\newtheorem{proposition}[theorem]{Proposition}

\def \N {{\mathbb N}}
\def \Z {{\mathbb Z}}
\def \C {{\mathbb C}}
\def \Q {{\mathbb Q}}

\def \P {{\mathbb P}}

\def \e {{\epsilon}}
\begin{document}

\title{Multiple Dedekind Zeta Values\\  are\\ Periods of Mixed Tate Motives}
\author{Ivan Horozov}
\maketitle

\begin{abstract} Recently, the author defined multiple Dedekind zeta values \cite{MDZF} associated to a number $K$ field and a cone $C$. These objects are number theoretic analogues of multiple zeta values.
In this paper we prove that every multiple Dedekind zeta value over any number field $K$ is a period of a mixed Tate motive. Moreover, if $K$ is a totally real number field, then we can choose a cone $C$ so that every multiple Dedekind zeta associated to the pair $(K;C)$ is unramified over the ring of algebraic integers in $K$. In \cite{period}, the author proves similar statements in the special case of a real quadratic fields for a particular type of a multiple Dedekind zeta values.

The mixed motives are defined over $K$ in terms of a the Deligne-Mumford compactification of the moduli space of curves of genus zero with $n$ marked points.
 
\end{abstract}




\section{Introduction}
The Riemann zeta function 
\[\zeta(s)=\sum_{n>0}\frac{1}{n^s}\]
is widely used in number theory, algebraic geometry and quantum field theory.
Euler's multiple zeta values
\[\zeta(s_1,\dots,s_m)=\sum_{0<n_1<\dots<n_m}\frac{1}{n_1^{s_1}\dots n_m^{s_m}},\]
where $s_1,\dots,s_m$ are positive integers and $s_m\geq 2$,
appear  as values of some Feynman amplitudes, and in algebraic geometry, as periods of mixed Tate motives over $Spec(\Z)$ (see \cite{GM}, \cite{DG}, \cite{B1}, \cite{KZ}).

Dedekind zeta values
\[\zeta_K(s)=\sum_{\mathfrak{a}\neq(0)}\frac{1}{N(\mathfrak{a})^s},\]
are a generalization of the Riemann zeta function to a number field $K$.
In some Feynman amplitudes one of the summands is $\log(1+\sqrt{2})$ or $\log\left(\frac{1+\sqrt{5}}{2}\right)$. These values are essentially the residues at $s=1$ of Dedekind zeta functions over $\Q(\sqrt{2})$ and over $\Q(\sqrt{5})$, respectively. 
For $s=2,3,4,\dots$ the values $\zeta_K(s)$ are periods of mixed Tate motives over the ring of algebraic integers in $K$ with ramification only at the discriminant of $K$ (see \cite{B2}).

In \cite{MDZF}, the author has constructed multiple Dedekind zeta values, which are a generalization of Euler's multiple zeta values to number fields in the same way as Dedekind zeta values generalizes Riemann zeta values. For a quadratic number field $K$, the key examples of multiple Dedekind zeta values are
\begin{equation}
\label{eq mdzv}
\zeta_{K;C}(s_1,\dots,s_m)=
\sum_{\beta_1,\dots,\beta_m\in C}
\frac{1}{N(\beta_1)^{s_1}N(\beta_1+\beta_2)^{s_2}\cdots N(\beta_1+\dots+\beta_m)^{s_m}},
\end{equation}
where $s_1,\dots,s_m$ are positive integers and $s_m\geq 2$ and
$C$ is a cone generated by totally positive algebraic integers $e_1,\dots,e_n$ in $K$  defined by
\[C=\N\{e_1,\dots,e_n\}=\{\gamma\in K\,|\, \gamma=a_1e_1+\dots a_ie_i,\mbox{ for positive integers } a_i\}.\]
Similar types of cones were considered by Zagier in \cite{Z1} and \cite{Z2}.

In \cite{MDZF}, the author has proven that multiple Dedekind zeta values can be interpolated to multiple Dedekind zeta functions, which have meromorphic continuation to all complex values of the variables $s_1,\dots,s_m$.

In this paper we prove the following two theorems.
\begin{theorem}
\label{thm}
Every multiple Dedekind zeta over any number field $K$ is a period of a mixed Tate motive over $K$. 
\end{theorem}

\begin{theorem}
\label{thm2}
If $K$ is a totally real field, then we can find a cone $C$ such that every multiple Dedekind zeta $\zeta_{K;C}(s_1,\dots,s_m)$ is a period of mixed Tate motive, which is unramified over any prime.
\end{theorem}
\section{Background}

\subsection{Multiple zeta values and iterated integrals}
The Riemann zeta function at the value $s=2$ can be expressed in term of an iterated integral in the following way

\begin{multline*}
\int_0^1\left(\int_0^y\frac{dx}{1-x}\right)\frac{dy}{y}
=\int_0^1\left(\int_0^y(1+x+x^2+x^3\dots)dx\right)\frac{dy}{y}\\
=\int_0^1\left(y+\frac{y^2}{2}+\frac{y^3}{3}+\frac{y^4}{4}+\dots \right)\frac{dy}{y}
=\left. y+\frac{y^2}{2^2}+\frac{y^3}{3^2}+\frac{y^4}{4^2}\cdots\right|_{y=0}^{y=1}\\
=1+\frac{1}{2^2}+\frac{1}{3^2}+\frac{1}{4^2}\dots
=\zeta(2).
\end{multline*}

Let us examine the domain of integration of the iterated integral. Note that $0<x<y$ and $0<y<1$. We can put both inequalities together. Then we obtain the domain $0<x<y<1$, which is a simplex. Thus, we can express the iterated integral as
\[
\zeta(2)=\int_0^1\left(\int_0^y\frac{dx}{1-x}\right)\frac{dy}{y}=\int_{0<x<y<1}\frac{dx}{1-x}\wedge\frac{dy}{y}.
\]

Moreover, Goncharov and Manin \cite{GM} have expressed all multiple zeta values of weight $m$ as periods of motives  related to the moduli space of curves of genus zero with $m+3$ marked points, $\mathcal{M}_{0,m+3}$. In particular, $\zeta(2)$ can be expressed as a period of the motive $H^2(\overline{\mathcal{M}}_{0,5}-A,B-A\cap B)$ by pairing of $[\Omega_A]\in Gr^W_4H^2(\overline{\mathcal{M}}_{0,5}-A)$ for $\Omega_A=\frac{dx}{1-x}\wedge\frac{dy}{y}$, 
with $[\Delta_B]\in \left(Gr^W_0H^2(\overline{\mathcal{M}}_{0,5}-B)\right)^\vee$. The Deligne-Mumford compactification $\overline{\mathcal{M}}_{0,5}$ of the moduli space $\mathcal{M}_{0,5}$ can be obtained by three blow-ups of $\P^1\times\P^1$ at the points $(0,0)$, $(1,1)$ and $(\infty,\infty)$. 
Let us name the exceptional divisors at the three points by $E_0$, $E_1$ and $E_\infty$, respectively. Then
$A=(x=1)\cup(y=0)\cup (x=\infty)\cup (y=\infty)\cup E_\infty$ and $B=(x=0)\cup(x=y)\cup(y=1)\cup E_0\cup E_1$.

Similarly, one can express $\zeta(3)$ and $\zeta(1,2)$ as iterated integrals
\begin{eqnarray*}
\zeta(3)&=&\int_0^1\left(\int_0^z\left(\int_0^y\frac{dx}{1-x}\right)\frac{dy}{y}\right)\frac{dz}{z}
=\int_{0<x<y<z<1}\frac{dx}{1-x}\wedge\frac{dy}{y}\wedge\frac{dz}{z},\\
\\
\zeta(1,2)&=&\int_0^1\left(\int_0^z\left(\int_0^y\frac{dx}{1-x}\right)\frac{dy}{1-y}\right)\frac{dz}{z}
=
\int_{0<x<y<z<1}\frac{dx}{1-x}\wedge\frac{dy}{1-y}\wedge\frac{dz}{z}.
\end{eqnarray*}
Again, $\zeta(3)$ and $\zeta(1,2)$ can be expressed as periods of motives related to $\mathcal{M}_{0,6}$.
In the same paper, Goncharov and Manin prove that the motives associated to multiple zeta values (MZVs) are mixed Tate motives unramified over $Spec(\Z)$.

A few years later, Francis Brown \cite{B1} proved that periods of mixed Tate motives unramified over $Spec(\Z)$ can be expressed as a $\Q$-linear combination of MZVs times an integer power of $2\pi i$.


\subsection{Multiple Dedekind zeta values (MDZVs) and iterated integrals on membranes}

We recall the construction of MDZVs. 
Let ${\mathcal{O}}_K$ be the ring of integers in a number field $K$ of degree $n$ over $\Q$.

Denote by $\sigma_1,\dots, \sigma_n$ all the embedding of $K$ into the complex numbers $\C$.
And let  $e_1,\dots, e_n$  be elements of  ${\mathcal O}_K$ such that
\begin{enumerate}
\item $e_i\in {\mathcal{O}}_K$ for all $i$
\item $(e_1,\dots,e_n)$ forms a basis of $K$ over $\Q$
\item $Re(\sigma_j(e_i)\geq 0$ for all $i$ and $j$
\end{enumerate}
Let $C$ be the cone defined as $\N$-linear combinations of $e_1,\dots,e_2$, that is,
\[C=\{\alpha\in {\mathcal{O}}_K\,|\, \gamma=a_1e_1+\cdots+a_ne_n,\mbox{ for }a_1,\dots,a_n\in\N\}.\]
Let \[f_0(C;t_1,\dots,t_n)=\sum_{\alpha\in C}\exp\left[-\sum_{i=1}^n t_j\sigma_j(\alpha)\right].\]
We express $\zeta_{K;C}(2)$, $\zeta_{K;C}(3)$ and $\zeta_{K;C}(1,2)$  as iterated integrals on a membrane.
See \cite{MDZF} and \cite{Hilbert}, for more examples and properties of  these constructions.
We have
\begin{eqnarray}
\label{int z2}
&&
\int_0^\infty\cdots\int_0^\infty\left(\int_{u_1}^\infty\cdots \int_{u_n}^\infty f_0(C;t_1,\dots,t_n)dt_1\dots dt_n\right)du_1\dots du_n\\
\nonumber
&&
=\int_0^\infty\cdots\int_0^\infty
\left(
\int_{u_1}^\infty\cdots \int_{u_n}^\infty 
\sum_{\alpha\in C}
\exp
\left[
-\sum_{j=1}^n t_j\sigma_j(\alpha) 
\right]
dt_1\dots dt_n
\right)
du_1\dots du_n\\
\nonumber
&&
=\int_0^\infty\cdots\int_0^\infty \left(\sum_{\alpha\in C}\left(\prod_{j=1}^n\exp\left[- t_j\sigma_j(\alpha)\right]\right)dt_1\dots dt_n\right)
du_1\dots du_n\\
\nonumber
&&
=\int_0^\infty\cdots\int_0^\infty \sum_{\alpha\in C}\prod_{j=1}^n\frac{\exp\left[- u_j\sigma_j(\alpha)\right]}{\prod_j\sigma_j(\alpha)}du_1\dots du_n\\
\nonumber
&&=\sum_{\alpha\in C}\frac{1}{(\prod_j\sigma_j(\alpha))^2}
=\sum_{\gamma\in C}\frac{1}{N(\alpha)^2}\\
\nonumber
&&=\zeta_{K;C}(2).
\end{eqnarray}

Similarly,
\begin{eqnarray*}
&&\int_0^\infty\cdots\int_0^\infty
\left(\int_{v_1}^\infty\cdots\int_{v_n}^\infty
\left(\int_{u_1}^\infty\cdots\int_{u_n}^\infty\right.\right.\\
&& 
f_0(C;t_1,\dots,t_n)dt_1\cdots dt_n\bigg)
du_1\cdots du_n\bigg)
dv_1\cdots dv_n\\
&&=\int_0^\infty\cdots\int_0^\infty
\left(
\int_{v_1}^\infty\cdots\int_{v_n}^\infty
\left(
\int_{u_1}^\infty\cdots \int_{u_n}^\infty
\right.
\right.\\ 
&&
\left.
\left.
\sum_{\alpha\in C}
\exp
\left[
-\sum_{j=1}^n t_j\sigma_j(\alpha) 
\right]
dt_1\dots dt_n
\right)
du_1\dots du_n
\right)
dv_1\cdots dv_n\\
&&=\sum_{\alpha\in C}
\int_0^\infty\cdots\int_0^\infty
\left(
\int_{v_1}^\infty\cdots\int_{v_n}^\infty
\prod_{j=1}^n\frac{\exp
\left[
- u_j\sigma_j(\alpha)
\right]}{\sigma_j(\alpha)}
du_1\dots du_n
\right)
dv_1\dots dv_n
\\
\\
&&=\sum_{\alpha\in C}
\int_0^\infty\cdots\int_0^\infty
\prod_{j=1}^n\frac{\exp
\left[
- v_j\sigma_j(\alpha)
\right]}{\sigma_j(\alpha)^2}
dv_1\dots dv_n
=\sum_{\alpha\in C}\frac{1}{N(\alpha)^3}\\
&&=\zeta_{K;C}(3).
\end{eqnarray*}

\newpage
We recall the simplest type of multiple Dedekind zeta value

\begin{eqnarray*}
&&\int_0^\infty\cdots\int_0^\infty
\left(\int_{v_1}^\infty\cdots\int_{v_n}^\infty
\left(\int_{u_1}^\infty\cdots\int_{u_n}^\infty\right.\right.\\
&& 
f_0(C;t_1,\dots,t_n)dt_1\cdots dt_n\bigg)
f_0(C;u_1,\dots,u_n)du_1\cdots du_n\bigg)
dv_1\cdots dv_n\\
&&=\int_0^\infty\cdots\int_0^\infty
\left(
\int_{v_1}^\infty\cdots\int_{v_n}^\infty
\left(
\int_{u_1}^\infty\cdots \int_{u_n}^\infty
\right.
\right.\\ 
&&
\left.
\left.
\sum_{\alpha\in C}
\exp
\left[
-\sum_{j=1}^n t_j\sigma_j(\alpha) 
\right]
dt_1\dots dt_n
\right)
\sum_{\beta\in C}
\exp
\left[
-\sum_{j=1}^n u_j\sigma_j(\beta) 
\right]
du_1\dots du_n
\right)
dv_1\cdots dv_n\\
&&=\sum_{\alpha,\beta\in C}
\int_0^\infty\cdots\int_0^\infty
\left(
\int_{v_1}^\infty\cdots\int_{v_n}^\infty
\prod_{j=1}^n\frac{\exp
\left[
- u_j\sigma_j(\alpha)
\right]}{\sigma_j(\alpha)}
\exp
\left[
- u_j\sigma_j(\beta)
\right] du_1\dots du_n
\right)
dv_1\dots dv_n
\\
&&=\sum_{\alpha,\beta\in C}
\int_0^\infty\cdots\int_0^\infty
\left(
\int_{v_1}^\infty\cdots\int_{v_n}^\infty
\prod_{j=1}^n\frac{\exp
\left[
- u_j\sigma_j(\alpha+\beta)
\right]}{\sigma_j(\alpha)}
du_1\dots du_n
\right)
dv_1\dots dv_n
\\
&&=\sum_{\alpha,\beta\in C}
\int_0^\infty\cdots\int_0^\infty
\prod_{j=1}^n\frac{\exp
\left[
- v_j\sigma_j(\alpha+\beta)
\right]}{\sigma_j(\alpha)\sigma_j(\alpha+\beta)}
dv_1\dots dv_n
=\sum_{\alpha\in C}\frac{1}{N(\alpha)^1N(\alpha+\beta)^2}\\
&&=\zeta_{K;C}(1,2).
\end{eqnarray*}
\section{Transition to Algebraic Geometry}
We may write the infinite sum in the definition of $f_0$ as a product of $n$ geometric series as follows.

\begin{lemma}
Let $x_i=e^{-t_i}$ for $i=1,2,\dots,n$.
Then 
$e^{-t_j\sigma_j(e_i)}=x_j^{\sigma_j(e_i)}$ and
\begin{equation}
\label{alg}
f_0(C;t_1,\dots, t_n)
=
\prod_{i=1}^n
\left(\frac{\prod_{j=1}^nx_j^{\sigma_j(e_i)}}{1-\prod_{j=1}^nx_j^{\sigma_j(e_i)}}\right)
\end{equation}
\end{lemma}
{\bf Proof.} We simplify the function $f_0$ by expressing it in terms of products:
\begin{eqnarray*}
&&f_0(C;t_1,\dots,t_n)
=
\sum_{\alpha\in C}
\exp\left[-\sum_{j=1}^n\sigma_j(\alpha)t_j\right]\\
&&=
\sum_{a_1=1}^\infty
\dots
\sum_{a_n=1}^\infty 
\exp\left[-\sum_{j=1}^n t_j[a_1\sigma_j(e_1)+\cdots+a_n\sigma_j(e_n)]\right]\\
&&=
\sum_{a_1=1}^\infty
\dots
\sum_{a_n=1}^\infty 
\exp[-a_1[t_1\sigma_1(e_1)+\cdots+ t_n\sigma_n(e_1)]]\times\\
&&\cdots\times
\exp[-a_n[t_1\sigma_1(e_n)+\cdots+ t_n\sigma_n(e_n)]]\\
&&=
\frac{\exp\left[-[t_1\sigma_1(e_1)+\cdots+ t_n\sigma_n(e_1)]\right]}{1-\exp[-[t_1\sigma_1(e_1)+\cdots+ t_n\sigma_n(e_1)]]}\times\\
&&\cdots\times
\frac{\exp[-[t_1\sigma_1(e_n)+\cdots+ t_n\sigma_n(e_n)]]}{1-\exp[-[t_1\sigma_1(e_n)+\cdots+ t_n\sigma_n(e_n)]]}=\\
&&=
\prod_{i=1}^n
\frac{\exp[-\sum_{j=1}^nt_j\sigma_j(e_i)]}{1-\exp\left[-\sum_{j=1}^nt_j\sigma_j(e_i)\right]}\\
&&=
\prod_{i=1}^n
\frac{\prod_{j=1}^n\exp[-t_j\sigma_j(e_i)]}{1-\prod_{j=1}^n\exp[-t_j\sigma_j(e_i)]}.\,\,\,\,\,\,\square
\end{eqnarray*}

\subsection{The Algebraic Exponent}
We are going to define new variables $x_{ij}$, so that when we express $f_0(C;t_1,\dots,t_n)$ in terms of $x_{ij}$, then we obtain a rational function. Intuitively $x_{ij}=x_j^{\sigma_j(e_i)}$. 
To achieve that, we need to define algebraically $f^\gamma(x)=x^\gamma$ where $\gamma\in {\cal O}_K$.
We follow similar ideas as in the announcement \cite{period}.
Let ${\cal O}_K=\Z\{\mu_1,\mu_2,\dots,\mu_n\}$ as a $\Z$-module, where $\mu_1=1$.
Let $(c_{ij})$ by the $n\times n$-matrix associated to $\gamma$ in the basis $(\mu_1,\dots,\mu_n)$, that is,
\[\gamma\mu_i=\sum_{k=1}^n c_{ik}\mu_k.\]
We define a function $f^{\gamma}$ corresponding to raising to a power $\gamma$ by sending monomials in the variables $y_1,\dots,y_n$ to monomials in the same variables. Let
\[f^{\gamma}(y_{i})=\prod_{k=1}^n y_{k}^{c_{ik}}.\]
\begin{lemma} Iterated application of the above definition of exponentiation has the following property:
\[f^{\beta}f^{\alpha}(y_{i})=f^{\alpha\beta}(y_{i}).\]
\end{lemma}
{\bf Proof:}
It follows from the fact that $\gamma\mapsto (c_{ik})$ is a representation of the ring ${\cal O}_K$ as an endomorphism of $\Z^n$.
More precisely, let $\alpha \mapsto (a_{ij})$, $\beta\mapsto (b_{jk})$ and  $\gamma\mapsto (c_{ik})$. If $\alpha\beta=\gamma$
then $\sum_j a_{ij}b_{jk}=c_{ik}$. Thus,
\begin{eqnarray*}
f^{\beta}f^{\alpha}(y_{i})=f^\beta\left(\prod_j y_j^{a_{ij}}\right)
=\prod_{j,k} y_{k}^{a_{ij}b_{jk}}=
\prod_{k} y_{k}^{\sum_ja_{ij}b_{jk}}=
\prod_{j,k} y_{k}^{c_{ik}}=f^{\alpha\beta}(y_{i}).\,\,\,\square
\end{eqnarray*}

Now that we have defined an algebraic power of a variable, we return to expressing $f_0$ as a rational function.
Let \[x_{ij}=f^{\sigma(e_i)}(x_j).\]
Intuitively, $x_{ij}=x_j^{\sigma_j(e_i)}$. 

Then\[
f_0(C;t_1,\dots,t_n)=\prod_{i=1}^n
\left(\frac{\prod_{j=1}^nx_j^{\sigma_j(e_i)}}{1-\prod_{j=1}^nx_j^{\sigma_j(e_i)}}\right)\]
can be written formally as
\begin{equation}
\label{eq f0 alg}
f_0(C;t_1,\dots,t_n)=
\prod_{i=1}^n
\left(\frac{\prod_{j=1}^n f^{\sigma_j(e_i)}(x_{j})}{1-\prod_{j=1}^n f^{\sigma_j(e_i)}(x_{j})}\right)=
\prod_{i=1}^n
\left(\frac{\prod_{j=1}^n x_{ij}}{1-\prod_{j=1}^n x_{ij}}\right).
\end{equation}


\section{Multiple Dedekind Zeta Values, Differential Forms and Rational Functions}
Let us recall the definition of a multiple Dedekind zeta value (see \cite{MDZF}).
Let 
\begin{eqnarray*}
&&\alpha_0(t_1,\dots,t_n)=dt_1\wedge\cdots\wedge dt_n\\
&&\alpha_1(t_1,\dots,t_n)=f_0(t_1,\dots,t_n)dt_1\wedge\cdots\wedge dt_n.
\end{eqnarray*}
The definition of a multiple Dedekind zeta value \eqref{eq mdzv} is as follows.
\begin{equation}
\label{def int}
\zeta_{K;C}(s_1,\dots,s_m)=\int_{\Delta}\bigwedge_{k=1}^M\alpha_{\e_k}(t_{1,k},\cdots,t_{n,k}),
\end{equation}
where 
\begin{enumerate}
\item $M=s_1+\cdots+s_d$; 
\item$\Delta=\Delta_1\times\cdots\times \Delta_n$, is an $n$-fold product of $M$-simplices $\Delta_1,\dots,\Delta_n$
\item $\Delta_j$ is a simplex consisting of points $(t_{j,1},t_{j,2},\cdots,t_{j,M})$ such that $t_{j,1}>t_{j,2}>\cdots >t_{j,M}>0$;
\item
the indecies $\e_k$ have values $0$ or $1$ and 
\begin{eqnarray*}
&\e_1=1, &\e_1=\dots=\e_{s_1}=0\\
&\e_{s_1+1}=1, &\e_{s_1+2}=\dots=\e_{s_1+s_2}=0\\
&\cdots & \\
&\e_{s_1+\dots+s_{d-1}+1}=1, &\e_{s_1+\dots+s_{d-1}+2}=\dots=\e_{s_1+\dots+s_{d}}=0
\end{eqnarray*}
\end{enumerate}

{\bf Definition.}
Let $z_i=\prod_{j=1}^n x_{ij}$.
We define the differential forms $\omega_0$ and $\omega_1$ by  
\begin{eqnarray}
\label{eq alg diff forms}
&&\omega_0(z_1,\dots,z_n)=\bigwedge_{i=1}^n\frac{dz_i}{z_i}\\
&&\omega_1(z_1,\dots,z_n)=\bigwedge_{i=1}^n \frac{dz_i}{1-z_i}
\end{eqnarray}
\begin{proposition}
Evaluate $x_{ij}$ at $e^{-t_j^{\sigma{e_i}}}$.   Then
\begin{eqnarray*}&&\omega_0(z_1,\dots,z_n)=\sqrt{\Delta}\cdot\alpha_0(t_1,\dots,t_n)\\
&&\omega_1(z_1,\dots,z_n)=\sqrt{D}\cdot\alpha_1(t_1,\dots,t_n),
\end{eqnarray*}
where $\sqrt{D}=\det(\sigma_j(e_i))$ and $D$ is an integer multiple of the discriminant.
\end{proposition}
{\bf Proof.}
If we evaluate $x_{ij}$ at $e^{-t_j\sigma_j(e_i)}$, then
$z_i=\prod_j x_{ij}=\prod_je^{-t_j\sigma_j(e_i)}=e^{-\sum_j t_j\sigma_j(e_i)}$.
Then 
\begin{eqnarray}
\label{eq omega0}
\omega_0(z_1,\dots,z_n)&=&\bigwedge_i \frac{dz_i}{z_i}=det(\sigma_j(e_i))\bigwedge dt_i=\\
&=&
\nonumber
\sqrt{D}\cdot\alpha_0(t_1,\dots,t_n).
\end{eqnarray}
In that case, we also have
\[f_0(C;t_1,\dots,t_n)=\prod_i\frac{z_i}{1-z_i}.\]
Therefore,
\begin{eqnarray}
\label{eq omega1}
\omega_1(z_1,\dots,z_n)&=&\bigwedge_i \frac{dz_i}{1-z_i}=\left(\prod_i\frac{z_i}{1-z_i}\right)\bigwedge_i \frac{dz_i}{z_i}=\\
&=&
\nonumber
\det(\sigma_j(e_i))f_0(t_1,\dots,t_n)\bigwedge dt_i=\\
&=&
\nonumber 
\sqrt{D}\cdot\alpha_1(t_1,\dots,t_n).
\end{eqnarray}

\begin{proposition}
We have the following relation between multiple Dedekind zeta values of the differential forms $\omega_0$ and $\omega_1$
\[\int_{\Delta}\bigwedge_{k=1}^M\omega_{\e_k}(t_{1,k},\cdots,t_{n,k})=\left(\sqrt{D}\right)^M\zeta_{K;C}(s_1,\dots,s_m).\]
\end{proposition}
{\bf Proof.}
It follows directly from Equations \eqref{def int}, \eqref{eq omega0} and \eqref{eq omega1}.
$\,\,\,\square$
\section{Tangential base points}

Let $y=e^{-bt }$ and $z=e^{-ct}$, where $b$ and $c$ are complex numbers such that $Re(b)>0$ and $Re(c)>0$, and $|b|\neq |c|.$
Then
\[
\lim_{t\rightarrow +\infty}\frac{dy}{dz}=
\lim_{t\rightarrow +\infty}\frac{de^{-bt}}{de^{-ct}}=
\lim_{t\rightarrow +\infty}\frac{de^{ct}}{de^{bt}}=q,
\]
where 
\[p=\left\{
\begin{tabular}{rll}
$+\infty$ & or $[0:1]$ & if  $b<c$\\
$0$  & or $[0:1]$ &if  $b>c$
\end{tabular}
\right.
.\]

Also
\[
\lim_{t\rightarrow 0}\frac{dy}{dz}=
\lim_{t\rightarrow 0}\frac{de^{-bt}}{de^{-ct}}=
\lim_{t\rightarrow 0}\frac{-b e^{-bt}}{-c e^{-ct}}=
\frac{b}{c}
\]

Let $\gamma:[0,1]\rightarrow {\cal M}_{5}$,
by sending $t$ to $(y,z)$, where $y=e^{-bt }$ and $z=e^{-ct}$.
For a vector $v=(a,b)$, consider $[v]=[a:b]$ as an element of $\P^1$.

We have proven the following lemma.
\begin{lemma}

(a)  \[\lim_{t \rightarrow \infty}\left[\frac{d\gamma}{dt}\right]
=\left\{
\begin{tabular}{rll}
$[0:1]$ & if  $b<c$\\
$[1:0]$  & if $b>c$
\end{tabular}
\right.
\]
Moreover, the limit is well defined for for any distinct positive real numbers $b$ and $c$.
  
(b) \[\lim_{t \rightarrow 0}\left[\frac{d\gamma}{dt}\right]=[b:c].\]
\end{lemma}

Let $t=t_j$. Let also $b=\sigma_j(e_i)$, $c=\sigma_j(e_k)$.
Then \[\lim_{t_j\rightarrow 0}\left[\frac{dx_{i_1,j}}{dx_{i_2,j}}\right]=[\sigma_j(e_{i}):\sigma_j(e_{k})].\] 
We define \[[q(i,k)]=[\sigma_j(e_{i}):\sigma_j(e_{k})].\]
And
\[\lim_{t_j\rightarrow \infty}\left[\frac{dx_{i,j}}{dx_{k,j}}\right]
=
\left\{
\begin{tabular}{rll}
$[0:1]$ & if  $\sigma_j(e_i)<\sigma_j(e_k)$\\
$[1:0]$  & if $\sigma_j(e_i)>\sigma_j(e_k)$
\end{tabular}
\right.
\]
Let
\[
[p(i,k)]=\left\{
\begin{tabular}{rll}
$[0:1]$ & if  $\sigma_j(e_i)<\sigma_j(e_k)$\\
$[1:0]$  & if $\sigma_j(e_i)>\sigma_j(e_k)$
\end{tabular}
\right.
\]
More generally, let 
$[p(i_0,\dots,i_r)]=[a_0:\cdots:a_r]$ be a point on $\P^r(\Q)$ with $a_i$ being $0$ or $1$, where all $a_i$'s are zero except the one whose index $c$, for $a_c$,  is such that $|\sigma_j(e_c)|$ is a maximal value among the elements $|\sigma_j(e_0)|,\dots,|\sigma_j(e_r)|$.
\section{Multiple Dedekind Zeta Values and the Moduli Space ${\cal M}_{0,m\cdot n^2+3}$}
Let $z_{ik}=\prod_{j=1}^k x_{ij}.$ there are $n^2$ such variables.
If we consider multiple Dedekind zeta value of depth $m$ then we need $m\cdot n^2$ variables
$(z_{ikd})_{i,k,d=1,1,1}^{n,n,m}\in {\cal M}_{0,m\cdot n^2+3}$.
The dimension of the differential form $\Omega(A)$ is $mn$.
Let $\e_d$ be $0$ or $1$ for $d=1,2,\dots,m$. Let also $\e_1=1$ and $\e_m=0$.

\[\Omega(A)=\bigwedge_{d=1}^m \omega_{\e_d}.\]
where \[A=(z_{ind}=\e_{d})_{i,d=1,1}^{n,m}\cup (z_{ind}=\infty)_{i,d=1,1}^{n,m}\]
and 
$B=B_1\cap B_2\cap \cdot \cap B_n$, where $codim B_r=r$ contains the divisors
\[B_1=(z_{i,1,1}=0)_{i=1}^{n}\cup (z_{i,j,d}=z_{i,j,d+1})_{i,j,d=1,1,1}^{n,n,m-1} \cup (z_{i,1,m}=1)_{i=1}^{n}\]
 together with the intersection of boundary components of ${\overline{\cal M}}_{0,m\cdot n^2+3}-{\cal M}_{0,m\cdot n^2+3}$ containing the same variable or the same constant $0$ or $1$.  Besides codimension 1 components, $B$ also contains a codimension 2 components.

Let $B'_{i_1,i_2}$ be a quiasi-subvariety in the boundary of the Deligne-Mumford compactification that has coordinates with  \[[z_{i_1,1,1}:z_{i_2,1,1}]=[p(i_1,i_2)].\] in the blow-up of the intersection $(z_{i_1,j,1}=0)\cap (z_{i_2,j,1}=0)$.
Let $B''_{i_1,i_2}$ be a cycle in the boundary of the Deligne-Mumford compactification above the intersection 
$(z_{i,1,m}=1)\cap(z_{k,j,m}=1)$ such that the coordinated of $B'_{i_1,i_2}$ in the blowup are
\[[1-z_{i_1,1,1}:1-z_{1_2,1,1}]=[\sigma_j(e_{i_1}):\sigma_j(e_{i_2})].\]
The codimension $2$ components of $B$ are the union of all $B'_{i_1,i_2}$ and $B''_{i_1,i_2}$. That is
\[B_2=\bigcup_{i_1<i_2}\left(B'_{i_1,i_2}\cup B''_{i_1,i_2}\right) .\]





Now let us write $\omega_0(x_1,x_2)$ and $\omega_1(x_1,x_2)$, when we want to specify the dependence on the variables.
In fact, both forms depend also on $y_1$ and $y_2$; however, we will take care of that by choosing a region of integration together with tangential base points.
\begin{theorem} (a) Every multiple Dedekind zeta value over a field $K$ times a suitable multiple of a power of the discriminant of $K$ is a period of a mixed Tate motive over $K$. More precisely, 
\[\left(\sqrt{D}\right)^M\zeta_{K;C}(s_1,\dots,s_m)=\int_{\Delta}\bigwedge_{k=1}^M\omega_{\e_k}(t_{1,k},\cdots,t_{n,k})\]
is a period of \[H^{nM}(\overline{{\cal M}}_{0,n^2M+3}-A,B-(A\cap B)),\] 
$\overline{{\cal M}}_{0,n^2M+3}$ is the Delingne-Mumford compactification of the moduli space of curves of genus zero with $n^2M+3$ marked points, and
 $A$ and $B$ that consist of a union of lower dimensional moduli spaces of curves of genus zero with marked points.

(b) For any field $K$ of degree $n$ over $\Q$, 
with the property that $K$ has $n$ units $e_1,\dots,e_n$ 
that are linearly independent over $\Q$ and $|\sigma_j(e_i)\neq|\sigma_j(e_k)|$, 
we have the following stronger statement. 
In particular, when $K$ is a totally real number field, 
we choose a cone $C=\N\{e_1,\dots,e_n\}.$ 
Then for any positive integers $s_1,\dots,s_m$ with $s_m\geq 2$, we have that 
\[\left(\sqrt{D}\right)^M\zeta_{K;C}(s_1,\dots,s_m)\]
is a period of an unramified mixed Tate motive over the ring of algebraic integers ${\cal O}_K$ in the field $K$.
\end{theorem}

{\bf{Proof:}} In this proof we are going to follow closely the paper by Goncharov and Manin \cite{GM}.
The period will be a pairing between 
 $[\Omega_A]\in Gr^W_{2nM}H^{nM}(\overline{\mathcal{M}}_{0,n^2M+3}-A)$ and $[\Delta_B]\in \left(Gr^W_0H^{nM}(\overline{\mathcal{M}}_{0,n^2M+3}-B)\right)^\vee$ associated to a mixed Tate motive $H^{nM}(\overline{\mathcal{M}}_{0,n^2M+3}-A;B-A\cap B)$.
 
We have that $A$ and $B_1$ are defined over $\Z$. 
Moreover, any component and any intersection of components of $A$ and $B$ are 
isomorphic to a moduli space ${\cal M}_{0,N}$ for some $N$. 
The component $B_2,B_3,\dots$ are defined over the field $K$. 
Moreover, any intersection is isomorphic to  ${\cal M}_{0,N/_K}$ for some $N$. Thus all multiple Dedekind zeta values are mixed Tate motives over the field of definition $K$.

If $e_1,\dots,e_n$ are unit in ${\cal O}_K$, which are linearly independent over $\Q$, then all $[q(i_1,i_2,j)]$, $[q(i_1,i_2,i_3,j])$ etc., have coordinates $0$ or units. Then, the component $B_2,B_3,\dots$ are defined over the ring ${\cal O}_K$. 
Moreover, any intersection is isomorphic to  ${\cal M}_{0,N/_{{\cal O}_K}}$ for some $N$. 

We have that $H^i(\overline{\mathcal{M}}_{0,N})$ is a mixed Tate motive over $Spec(\Z)$. This implies that
$H^i(\overline{\mathcal{M}}_{0,N})$ is a mixed Tate motive over $Spec(\mathcal{O}_K)$. we obtain that the motivic cohomology of the components of $B$ are mixed Tate motives. Using Proposition 1.7 from Deligne and Goncharov, \cite{DG}, we conclude that for $l\neq char (\nu)$ the $l$-adic cohomology of the reduction of $B_j$  modulo $\nu$ of the motive $H^i(B_j)$ is unramified for any component $B_j$ of $B$ , since $B_j$ is isomorphic to $\overline{\mathcal{M}}_{0,N}$ over $Spec(\mathcal{O}_K)$  for some $N$. We conclude that for $l\neq char (\nu)$ the $l$-adic cohomology of the reduction modulo any $\nu\in Spec(\mathcal{O}_K)$ of the motive $H^{nM}(\overline{\mathcal{M}}_{0,n^2M+3}-A;B-A\cap B)$ is unramified. Thus,
$H^{nM}(\overline{\mathcal{M}}_{0,n^2M+3}-A;B-A\cap B)$ is a mixed Tate motive unramified over $Spec(\mathcal{O}_K)$.

\section*{Acknowledgements}
I would like to express great respect to Emma Previato for her work and style in research and in communication. Without her interest in my work, this paper would not have been completed.
I am very thankful to Yuri Manin and Don Zagier for conversations during the preparation this paper and for their interest.

Part of this work was finished at the Max Planck Institute for Mathematics and I am grateful for the excellent working conditions and the financial support provided. Support for this project was also provided by a PSC-CUNY Award, jointly funded by The Professional Staff Congress and The City University of New York.





City University of New York,
Bronx Community College
2155 University Avenue, Bronx,
New York 10453 , U.S.A.;
ivan.horozov@bcc.cuny.edu

\end{document}